\documentclass[a4paper,12pt,twoside,leqno,final]{amsart}
\usepackage{amsmath}
\usepackage{amssymb}

\newtheorem{thm}{Theorem}[section]
\newtheorem{lem}[thm]{Lemma}

\theoremstyle{definition}

\newcommand{\C}{{\mathbb C}}
\newcommand{\D}{{\mathbb D}}

\newcommand{\T}{{\mathbb T}}
\newcommand{\Z}{{\mathbb Z}}

\newcommand{\La}{\Lambda}

\newcommand{\eps}{\varepsilon}

\newcommand{\f}{\frac}
\newcommand{\ov}{\overline}
\newcommand{\al}{\alpha}

\newcommand{\ga}{\gamma}

\newcommand{\de}{\delta}

\newcommand{\ze}{\zeta}

\newcommand{\ph}{\varphi}

\DeclareMathOperator*{\esssup}{ess\,sup}
\numberwithin{equation}{section}

\title[Inner functions as strongly extreme points]
{Inner functions as strongly extreme points:\\ 
stability properties}
\author{Konstantin M. Dyakonov}
\address{Departament de Matem\`atiques i Inform\`atica, IMUB, BGSMath, Universitat de Barcelona, Gran Via 585, E-08007 Barcelona, Spain}
\address{ICREA, Pg. Llu\'is Companys 23, E-08010 Barcelona, Spain}
\email{konstantin.dyakonov@icrea.cat}
\keywords{Bounded analytic functions, inner functions, extreme points, strongly extreme points}
\subjclass[2010]{30H05, 30H10, 30J05, 46A55}
\thanks{Supported in part by grant PID2021-123405NB-I00 from El Ministerio de Ciencia e Innovaci\'on (Spain) and grant 2021 SGR 00087 from AGAUR (Generalitat de Catalunya).}

\begin{document}
\begin{abstract}
Given a Banach space $\mathcal X$, let $x$ be a point in $\text{\rm ball}(\mathcal X)$, the closed unit ball of $\mathcal X$. One says that $x$ is a {\it strongly extreme point} of $\text{\rm ball}(\mathcal X)$ if it has the following property: for every $\varepsilon>0$ there is $\delta>0$ such that the inequalities $\|x\pm y\|<1+\delta$ imply, for $y\in\mathcal X$, that $\|y\|<\varepsilon$. We are concerned with certain subspaces of $H^\infty$, the space of bounded holomorphic functions on the disk, that arise upon imposing finitely many linear constraints and can be viewed as finite-dimensional perturbations of $H^\infty$. It is well known that the strongly extreme points of $\text{\rm ball}(H^\infty)$ are precisely the inner functions, while the (usual) extreme points of this ball are the unit-norm functions $f\in H^\infty$ with $\log(1-|f|)$ non-integrable over the circle. Here we show that similar characterizations remain valid for our perturbed $H^\infty$-type spaces. Also, we investigate to what extent a non-inner function can differ from a strongly extreme point.
\end{abstract}

\maketitle

\section{Introduction and statement of results}

We write $\T$ for the circle $\{\ze\in\C:|\ze|=1\}$ and $m$ for its Haar measure. We then consider the space $L^\infty=L^\infty(\T,m)$ of all essentially bounded complex-valued functions on $\T$, equipped with the essential supremum norm $\|\cdot\|_\infty$, along with its subspace $H^\infty$. By definition, an $L^\infty$ function is in $H^\infty$ if it agrees a.e. on $\T$ with the boundary trace of a bounded holomorphic function on the disk 
$$\D:=\{z\in\C:|z|<1\},$$
where the boundary values are understood in the sense of nontangential convergence almost everywhere. Equivalently, 
$$H^\infty=\{f\in L^\infty:\,\widehat f(k)=0\,\,\,\text{\rm for }\,k=-1,-2,\dots\},$$
where $\widehat f(k)$ is the $k$th {\it Fourier coefficient} of $f$ given by 
$$\widehat f(k):=\int_\T\ov\ze^kf(\ze)\,dm(\ze),\qquad k\in\Z.$$
The norm on $H^\infty$, as well as on its subspaces that appear below, is always taken to be $\|\cdot\|_\infty$. See any of \cite{G, Hof, K} for a systematic treatment of $H^\infty$ and, more generally, of the Hardy spaces $H^p$. 

\par We shall be dealing with certain subspaces of $H^\infty$ that can be thought of as finite-dimensional perturbations of the whole space. Each of these is annihilated by a fixed finite subset of the dual space $(H^\infty)^*$. No special knowledge about $(H^\infty)^*$ is required for our purposes (see, however, \cite[Chapter V]{G} for a kind of description of this dual space). Let us only mention the most obvious, and \lq\lq tangible," functionals in $(H^\infty)^*$ that are induced---under the natural pairing---by $L^1$ functions on $\T$, or rather by the corresponding cosets.

\par Now, given a nonempty set $\Phi\subset(H^\infty)^*$, we define
$$H^\infty_\Phi:=\{f\in H^\infty:\,\ph(f)=0\,\,\text{\rm for all }\ph\in\Phi\}.$$
When $\#\Phi<\infty$ (which is the case that interests us), it seems reasonable to expect that $H^\infty_\Phi$ will be fairly similar to $H^\infty$, in some sense or another, and our aim here is to furnish a couple of \lq\lq geometric" results to that effect. Specifically, we show that the structure of the unit ball of $H^\infty_\Phi$ inherits quite a bit from that of $H^\infty$, as far as certain types of boundary points are concerned. To proceed with precise formulations, we need to introduce the appropriate geometric concepts, and we do this now in the abstract Banach space setting.

\par Let $\mathcal X=(\mathcal X,\|\cdot\|)$ be a Banach space, and let
$$\text{\rm ball}(\mathcal X):=\{x\in\mathcal X:\,\|x\|\le1\}$$
be its closed unit ball. Recall that an element $x$ of $\text{\rm ball}(\mathcal X)$ is said to be an {\it extreme point} thereof if it is not writable in the form $x=\f12(u+v)$ with two distinct points $u,v\in\text{\rm ball}(\mathcal X)$. Equivalently, $x$ is an extreme point of $\text{\rm ball}(\mathcal X)$ if and only if the only vector $y\in\mathcal X$ satisfying 
$$\max\{\|x+y\|,\,\|x-y\|\}\le1$$
is $y=0$. 

\par Furthermore, a point $x\in\text{\rm ball}(\mathcal X)$ is called a {\it strongly extreme point} of the ball if for every $\eps>0$ there is $\de>0$ such that, whenever $y\in\mathcal X$ is a vector with
\begin{equation}\label{eqn:pmlessoneplde}
\max\{\|x+y\|,\,\|x-y\|\}<1+\de,
\end{equation}
we have
\begin{equation}\label{eqn:nylesseps}
\|y\|<\eps.
\end{equation}
A slightly different, but equivalent, definition of a strongly extreme point (for a generic convex set) appears in \cite{M}; the current version is borrowed from \cite{CT}. Of course, every strongly extreme point of $\text{\rm ball}(\mathcal X)$ is extreme, and every extreme point has norm $1$. 

\par Yet another piece of terminology will be needed. Namely, given a point $x\in\text{\rm ball}(\mathcal X)$ and a number $\eps>0$, we say that $\eps$ is {\it $\mathcal X$-admissible} (or just {\it admissible}) {\it for} $x$ if there is a $\de>0$ that makes the implication \eqref{eqn:pmlessoneplde}$\implies$\eqref{eqn:nylesseps} true for any $y\in\mathcal X$. Obviously, if $\eps$ is admissible for $x$, then so is every $\eps'$ with $\eps'>\eps$. Also, it is clear that $x$ is a strongly extreme point of $\text{\rm ball}(\mathcal X)$ if and only if every $\eps>0$ is admissible for $x$.

\par It follows directly from the definitions that if $\mathcal X_0$ is a subspace of $\mathcal X$, carrying the same norm, and if $x\in\mathcal X_0$ is a point which is extreme (or strongly extreme) for $\text{\rm ball}(\mathcal X)$, then it is also extreme (resp., strongly extreme) for $\text{\rm ball}(\mathcal X_0)$. Sometimes, the converse happens to be true, meaning that every (strongly) extreme point of the unit ball in the smaller space, $\mathcal X_0$, is (strongly) extreme in the bigger one, $\mathcal X$. If so, we can speak of a certain {\it stability phenomenon} that occurs when passing from $\mathcal X$ to $\mathcal X_0$, in the sense that no new (strongly) extreme points emerge. 

\par Our first theorem establishes such a stability property (for both types of points) when the two spaces are $H^\infty$ and $H^\infty_\Phi$, with a finite set $\Phi\subset(H^\infty)^*$. Before stating it, we recall that a unit-norm function $f\in H^\infty$ is an extreme point of $\text{\rm ball}(H^\infty)$ if and only if 
\begin{equation}\label{eqn:logintegraldiv}
\int_\T\log(1-|f|)\,dm=-\infty;
\end{equation}
see \cite[Section V]{dLR} or \cite[Chapter 9]{Hof}. Furthermore, a theorem of Cima and Thomson \cite{CT} tells us that a function $f$ is a strongly extreme point of $\text{\rm ball}(H^\infty)$ if and only if $f$ is an {\it inner function} (i.e., $f\in H^\infty$ and $|f|=1$ a.e. on $\T$); see also \cite{MvR} for an extension of this last result to general uniform algebras. We now claim that similar characterizations are valid for our $H^\infty_\Phi$ spaces.

\begin{thm}\label{thm:exstex} Let $\Phi$ be a finite subset of $(H^\infty)^*$. Assume further that $f\in H^\infty_\Phi$ and $\|f\|_\infty=1$. Then 
\par\noindent{\rm(a)} $f$ is an extreme point of $\text{\rm ball}(H^\infty_\Phi)$ if and only if \eqref{eqn:logintegraldiv} holds;
\par\noindent{\rm(b)} $f$ is a strongly extreme point of $\text{\rm ball}(H^\infty_\Phi)$ if and only if it is an inner function.
\end{thm}

\par We mention in passing that things become different in the setting of the Hardy space $H^1$, where the stability phenomenon breaks down even for subspaces of codimension $1$. This will be explained in Section 3 below.

\par Regarding part (a) of Theorem \ref{thm:exstex}, we note that a tiny special case was previously proved in \cite{DAA}. There, the functionals that constitute $\Phi$ were of the form 
$$h\mapsto\widehat h(k_j),\qquad h\in H^\infty,$$
for some positive integers $k_1,\dots,k_N$, so that $H^\infty_\Phi$ was the space of $H^\infty$ functions with prescribed holes in the spectrum. 

\par Here, we are chiefly concerned with strongly extreme points, and we now turn to another stability issue related to statement (b) in the theorem above. Namely, bearing that statement in mind, we feel tempted to ask (albeit somewhat vaguely) to what extent a {\it non-inner} unit-norm function in $H^\infty_\Phi$ may fall short of being strongly extreme for $\text{\rm ball}(H^\infty_\Phi)$. This seems to call for a suitable quantitative refinement of (b), so we proceed in search of such a result. 

\par Given a finite set $\Phi\subset(H^\infty)^*$ and a function $f\in H^\infty_\Phi$ with $\|f\|_\infty=1$, we put
$$\eps_\Phi(f):=\inf\{\eps>0:\,\eps\,\,\text{\rm is }
H^\infty_\Phi\text{\rm -admissible for }f\}.$$
Thus, $\eps_\Phi(f)=0$ if and only if $f$ is a strongly extreme point of $\text{\rm ball}(H^\infty_\Phi)$.

\par Next, for a unit-norm function $f\in H^\infty$ and a number $\eta\in(0,1)$, we consider the sublevel set
\begin{equation}\label{eqn:sublevel}
E_\eta(f):=\{\ze\in\T:\,|f(\ze)|\le\eta\},
\end{equation}
defined up to a set of measure $0$. It should be noted that $f$ is inner if and only if $m(E_\eta(f))=0$ for each $\eta\in(0,1)$. Consequently, for $f$ non-inner, we might say that $f$ is far from (resp., close to) being an inner function if the number 
$$\inf\{\eta\in(0,1):\,m(E_\eta(f))>0\}$$
is small (resp., close to $1$). 

\par Further, once a finite set $\Phi\subset(H^\infty)^*$ and numbers $\eta,\ga\in(0,1)$ are fixed, we define 
$$\mathcal M(\Phi,\eta,\ga):=\left\{f\in H^\infty_\Phi:\,
\|f\|_\infty=1,\,\,m(E_\eta(f))\ge\ga\right\},$$
and we observe that every non-inner unit-norm function in $H^\infty_\Phi$ belongs to (at least) one of these sets. Finally, given a positive integer $N$ and numbers $\eta,\ga$ as above, we introduce the quantity 
$$\eps^*(N,\eta,\ga):=\inf_{\#\Phi=N}\inf\{\eps_\Phi(f):\,f\in\mathcal M(\Phi,\eta,\ga)\},$$
where the outer infimum is taken over the sets $\Phi\subset(H^\infty)^*$ of cardinality $N$.

\par We are now in a position to state our quantitative result, which offers a variation on the theme of Theorem \ref{thm:exstex}, part (b). 

\begin{thm}\label{thm:quantstex} Suppose that $0<\eta<1$, $0<\ga<1$, and $N$ is a positive integer. Then
\begin{equation}\label{eqn:estepsstar}
(c\ga)^N(1-\eta)\le\eps^*(N,\eta,\ga)\le(1-\eta^2)^{1/2}
\end{equation}
with a numerical constant $c>0$.
\end{thm}

\par An inequality of Nazarov, which is employed below to verify the left-hand estimate in \eqref{eqn:estepsstar}, provides further information on the constant $c$; specifically, it shows that one can take $c=\pi/(16e)$. 
\par The proofs of our results are given in the next section. We then conclude the paper by making some closing remarks and discussing several questions that puzzle us.

\section{Proofs of Theorems \ref{thm:exstex} and \ref{thm:quantstex}}

\noindent{\it Proof of Theorem \ref{thm:exstex}.} (a) The \lq\lq if" part is obvious, since $H^\infty_\Phi\subset H^\infty$. Indeed, \eqref{eqn:logintegraldiv} implies that $f$ is an extreme point of $\text{\rm ball}(H^\infty)$ and hence also of the smaller set $\text{\rm ball}(H^\infty_\Phi)$. 
\par To prove the \lq\lq only if" part, suppose that \eqref{eqn:logintegraldiv} fails, i.e., 
\begin{equation}\label{eqn:logintegralconv}
\int_\T\log(1-|f|)\,dm>-\infty.
\end{equation}
Let $N:=\#\Phi$, so that 
$$\Phi=\{\ph_1,\dots,\ph_N\}$$
for some pairwise distinct functionals $\ph_j$ from $(H^\infty)^*$. Also, let $\mathcal P_N$ stand for the set of (complex) polynomials of degree not exceeding $N$. Furthermore, \eqref{eqn:logintegralconv} enables us to consider the {\it outer function}, say $G$, with modulus $1-|f|$; this is defined by 
$$G(z):=\exp\left\{\int_\T\f{\xi+z}{\xi-z}\log(1-|f(\xi)|)\,dm(\xi)\right\},\qquad z\in\D,$$
and then extended to almost all $\ze\in\T$ by putting $G(\ze):=\lim_{r\to1^-}G(r\ze)$. We have $G\in H^\infty$ and $|G|=1-|f|$ a.e. on $\T$. 
\par Next, we claim that there exists a non-null polynomial $p\in\mathcal P_N$ for which $Gp\in H^\infty_\Phi$. To check this, we associate with each vector 
$$\al=(\al_0,\al_1,\dots,\al_N)\in\C^{N+1}$$
the polynomial 
$$p_\al(z):=\sum_{j=0}^N\al_jz^j\in\mathcal P_N$$
and let $T:\C^{N+1}\to\C^N$ be the linear map given by 
$$T\al=\left(\ph_1(Gp_\al),\dots,\ph_N(Gp_\al)\right).$$
Because the rank of $T$ does not exceed $N$, it follows from the rank-nullity theorem (see, e.g., \cite[p.\,63]{Axl}) that $\text{\rm Ker}\,T$, the null-space of $T$, has dimension at least $1$. In particular, $\text{\rm Ker}\,T$ is nontrivial. Now, if $\al$ is some (any) vector in $\text{\rm Ker}\,T\setminus\{0\}$, then the corresponding polynomial $p=p_\al$ is non-null and satisfies 
$$\ph_1(Gp)=\dots=\ph_N(Gp)=0,$$
whence $Gp\in H^\infty_\Phi$. 
\par Normalizing this $p$ so that $\|p\|_\infty=1$, we get
$$|f\pm Gp|\le|f|+|G|\,|p|\le|f|+|G|=1$$
a.e. on $\T$. Consequently, $f+Gp$ and $f-Gp$ are two distinct elements of $\text{\rm ball}(H^\infty_\Phi)$, so their midpoint (which is $f$) is not an extreme point of that ball. 

\smallskip (b) Again, the \lq\lq if" part is immediate. Indeed, the Cima--Thomson result from \cite{CT} tells us that every inner function is a strongly extreme point of $\text{\rm ball}(H^\infty)$, and we only need to combine this with the fact that $H^\infty_\Phi\subset H^\infty$. 

\par To prove the converse, assume that $f$ fails to be inner. This means that for some $\eta\in(0,1)$ we have $m(E_\eta(f))>0$, where $E_\eta(f)$ is the sublevel set defined by \eqref{eqn:sublevel}. Below, we put $E:=E_\eta(f)$. Also, we keep writing $N:=\#\Phi$ and recall the notation $\mathcal P_N$ for the space of polynomials of degree at most $N$. 
\par Further, we observe that there is a constant $C_{E,N}$ making the estimate 
\begin{equation}\label{eqn:estimq}
\|q\|_\infty\le C_{E,N}\sup_{\ze\in E}|q(\ze)|
\end{equation}
true for all $q\in\mathcal P_N$. Indeed, the supremum on the right-hand side defines a norm on the (finite-dimensional) space $\mathcal P_N$, which must be equivalent to $\|\cdot\|_\infty$. Finally, we consider the positive number
\begin{equation}\label{eqn:defepszero}
\eps_0:=(1-\eta)/C_{E,N}
\end{equation}
and we are going to show that $\eps_0$ is not $H^\infty_\Phi$-admissible for $f$. This will imply that $f$ is not a strongly extreme point for $\text{\rm ball}(H^\infty)$. 

\par Our goal is therefore to produce, for any given $\de>0$, a function $g\in H^\infty_\Phi$ such that
\begin{equation}\label{eqn:maxleonede}
\max\{\|f+g\|_\infty,\,\|f-g\|_\infty\}<1+\de
\end{equation}
and $\|g\|_\infty\ge\eps_0$. To this end, we fix $\de>0$ and let $\mathcal G$ be the outer function with modulus 
$$u:=(1-\eta)\cdot\chi_E+\f\de2\cdot\chi_{\T\setminus E}$$
on $\T$. Thus, $\mathcal G$ is defined as (the boundary trace of) 
$$\mathcal G(z):=\exp\left\{\int_\T\f{\xi+z}{\xi-z}\log u(\xi)\,dm(\xi)\right\},\qquad z\in\D,$$
so that $\mathcal G\in H^\infty$ and $|\mathcal G|=u$ a.e. on $\T$. Arguing as in the proof of part (a) above, we see that there exists a polynomial $p\in\mathcal P_N$ with $\|p\|_\infty=1$ for which $\mathcal Gp\in H^\infty_\Phi$. Setting $g:=\mathcal Gp$, we have then 
\begin{equation}\label{eqn:estont}
|f\pm g|\le|f|+|g|\le|f|+|\mathcal G|\quad\text{\rm a.e. on }\T.
\end{equation}
Recalling that $|f|\le\eta$ a.e. on $E$, while $|\mathcal G|$ takes the value $1-\eta$ (resp., $\de/2$) a.e. on $E$ (resp., on $\T\setminus E$), we further obtain 
\begin{equation}\label{eqn:estone}
|f|+|\mathcal G|\le\eta+(1-\eta)=1\quad\text{\rm a.e. on }E
\end{equation}
and
\begin{equation}\label{eqn:estontminuse}
|f|+|\mathcal G|\le1+\f\de2\quad\text{\rm a.e. on }\T\setminus E.
\end{equation}
It now follows from \eqref{eqn:estont}, \eqref{eqn:estone} and \eqref{eqn:estontminuse} that 
$$\|f\pm g\|_\infty\le1+\f\de2,$$
which obviously implies \eqref{eqn:maxleonede}. 

\par On the other hand, an application of \eqref{eqn:estimq} with $q=p$ yields
$$\sup_{\ze\in E}|p(\ze)|\ge1/C_{E,N}.$$
Consequently,
\begin{equation*}
\begin{aligned}
\|g\|_\infty\ge\esssup_{\ze\in E}|g(\ze)|&=\esssup_{\ze\in E}|\mathcal G(\ze)p(\ze)|\\
&=(1-\eta)\sup_{\ze\in E}|p(\ze)|\ge(1-\eta)/C_{E,N}=\eps_0,
\end{aligned}
\end{equation*}
whence $\|g\|_\infty\ge\eps_0$ as desired.\qed

\medskip Our next proof makes use of an explicit expression for the constant $C_{E,N}$ in \eqref{eqn:estimq}. This is provided by Lemma \ref{lem:turnaz} below, which is a restricted version of the Tur\'an--Nazarov inequality; see \cite[Theorem 1.4]{Naz}. 

\begin{lem}\label{lem:turnaz} There is an absolute constant $c\in(0,1)$ with the following property: whenever $n$ is a positive integer and $q$ is a continuous function on $\T$ with 
\begin{equation}\label{eqn:cardspeclen}
\#\{k\in\Z:\,\widehat q(k)\ne0\}\le n,
\end{equation}
the estimate 
\begin{equation}\label{eqn:estimturnaz}
\|q\|_\infty\le\left(\f1{cm(\mathcal E)}\right)^{n-1}\sup_{\ze\in\mathcal E}|q(\ze)|
\end{equation}
holds for any measurable set $\mathcal E\subset\T$ with $m(\mathcal E)>0$.
\end{lem}

\par The functions $q$ satisfying \eqref{eqn:cardspeclen} are, of course, trigonometric polynomials with at most $n$ terms. In particular, \eqref{eqn:estimturnaz} is valid when $q$ is an algebraic polynomial of degree not exceeding $n-1$. Nazarov's inequality, as proved in \cite{Naz}, yields a sharper version of \eqref{eqn:estimturnaz} which has the \lq\lq Wiener norm" $\sum_k|\widehat q(k)|$ in place of $\|q\|_\infty$ on the left-hand side; in addition, it tells us that the constant $\pi/(16e)$ is eligible as $c$. The case where $\mathcal E$ is an arc was treated much earlier by Tur\'an; see \cite{T}. A statement and proof of his original result can also be found in \cite[Part Two, Chapter 4]{HJ}, along with a discussion of Nazarov's improvement.

\medskip\noindent{\it Proof of Theorem \ref{thm:quantstex}.} Let $\Phi$ be a subset of $(H^\infty)^*$ consisting of $N$ elements, and let $f\in\mathcal M(\Phi,\eta,\ga)$. Writing $E:=E_\eta(f)$ as before, we thus have $m(E)\ge\ga$. Now, we know from the preceding proof that the number $\eps_0$, given by \eqref{eqn:defepszero}, is not $H^\infty_\Phi$-admissible for $f$; the constant $C_{E,N}$ appearing in \eqref{eqn:defepszero} was only required to make \eqref{eqn:estimq} true for all $q\in\mathcal P_N$. In particular, one possible choice is
$$C_{E,N}=\left(\f1{cm(E)}\right)^N$$
with the appropriate value of $c$, as asserted by Lemma \ref{lem:turnaz}. With this expression plugged in, \eqref{eqn:defepszero} takes the form 
$$\eps_0=(cm(E))^N(1-\eta);$$
and since this $\eps_0$ is not $H^\infty_\Phi$-admissible for $f$, neither is the smaller (or equal) number $(c\ga)^N(1-\eta)$. It follows that 
$$\eps_\Phi(f)\ge(c\ga)^N(1-\eta).$$
Taking the infimum over all $f\in\mathcal M(\Phi,\eta,\ga)$, and then over the sets $\Phi\subset(H^\infty)^*$ of cardinality $N$, we arrive at the left-hand inequality in \eqref{eqn:estepsstar}.

\par To prove the remaining part of \eqref{eqn:estepsstar}, it suffices to exhibit a set $\Phi\subset(H^\infty)^*$ with $\#\Phi=N$ and a function $f\in\mathcal M(\Phi,\eta,\ga)$ for which 
\begin{equation}\label{eqn:epssqrt}
\eps_\Phi(f)\le(1-\eta^2)^{1/2}.
\end{equation}
To this end, we define the functionals $\ph_k\in(H^\infty)^*$ with $k=1,\dots,N$ by
$$\ph_k(h)=\widehat h(k-1),\qquad h\in H^\infty,$$
and then put 
$$\Phi=\{\ph_1,\dots,\ph_N\}.$$
Further, we fix some (any) measurable set $E\subset\T$ with $m(E)=\ga$, and let $\mathcal F$ be the outer function such that
\begin{equation}\label{eqn:modcalf}
|\mathcal F|=\eta\cdot\chi_E+\chi_{\T\setminus E}\quad\text{\rm a.e. on }\T.
\end{equation}
Finally, we consider the function $f$ given by 
$$f(z):=z^N\mathcal F(z).$$
It is clear that $f\in H^\infty$ and $\widehat f(k)=0$ for $0\le k\le N-1$, whence $f\in H^\infty_\Phi$. Also, because $|f|=|\mathcal F|$ on $\T$, \eqref{eqn:modcalf} shows that $\|f\|_\infty=1$ and $E_\eta(f)=E$; therefore, $f\in\mathcal M(\Phi,\eta,\ga)$. 

\par We now claim that every number $\eps$ satisfying 
\begin{equation}\label{eqn:epsbiggersqrt}
\eps>(1-\eta^2)^{1/2}
\end{equation}
is $H^\infty_\Phi$-admissible for $f$. Once this is verified, \eqref{eqn:epssqrt} will follow readily, and the proof will be complete. 

\par Assuming \eqref{eqn:epsbiggersqrt}, we can clearly find a number $\de>0$ such that 
$$2\de+\de^2<\eps^2-(1-\eta^2),$$
or equivalently, 
\begin{equation}\label{eqn:deps}
(1+\de)^2-\eta^2<\eps^2.
\end{equation}
Now suppose that for some $g\in H^\infty_\Phi$ we have
\begin{equation}\label{eqn:gondon}
\|f\pm g\|_\infty<1+\de.
\end{equation}
Using the fact that $|f|\ge\eta$ a.e. on $\T$, the parallelogram identity and \eqref{eqn:gondon}, we then obtain
$$\eta^2+|g|^2\le|f|^2+|g|^2=\f12\left(|f+g|^2+|f-g|^2\right)\le(1+\de)^2,$$
whence
$$|g|^2\le(1+\de)^2-\eta^2$$
a.e. on $\T$. Together with \eqref{eqn:deps}, this yields
$$\|g\|_\infty\le\left\{(1+\de)^2-\eta^2\right\}^{1/2}<\eps.$$
To summarize, we have checked that $\eps$ is indeed $H^\infty_\Phi$-admissible for $f$, as long as \eqref{eqn:epsbiggersqrt} holds. Condition \eqref{eqn:epssqrt} is thereby established, and we are done.\qed

\section{Concluding remarks and open questions}

(1) Recall that, given a Banach space $\mathcal X$, a point $x$ of $\text{\rm ball}(\mathcal X)$ is said to be {\it exposed} for the ball if there exists a functional $\phi\in\mathcal X^*$ of norm $1$ such that 
$$\{y\in\text{\rm ball}(\mathcal X):\,\phi(y)=1\}=\{x\}.$$ 
A theorem of Amar and Lederer (see \cite{AL}) states that a function $f\in H^\infty$ with $\|f\|_\infty=1$ is an exposed point of $\text{\rm ball}(H^\infty)$ if and only if 
\begin{equation}\label{eqn:exposmeas}
m\left(\{\ze\in\T:\,|f(\ze)|=1\}\right)>0.
\end{equation}
In light of Theorem \ref{thm:exstex}, it seems natural to ask---assuming again that $\Phi$ is a finite subset of $(H^\infty)^*$---whether the exposed points of $\text{\rm ball}(H^\infty_\Phi)$ are precisely the unit-norm functions $f\in H^\infty_\Phi$ that obey \eqref{eqn:exposmeas}.

\smallskip (2) With regard to inequalities \eqref{eqn:estepsstar} in Theorem \ref{thm:quantstex}, it would be nice to narrow the gap between the two bounds for $\eps^*(N,\eta,\ga)$ so as to make the estimates asymptotically sharp.

\smallskip (3) While Theorem \ref{thm:exstex} says that no new extreme (or strongly extreme) points emerge when passing from $H^\infty$ to $H^\infty_\Phi$, provided that $\#\Phi<\infty$, no such stability occurs in the context of the Hardy space $H^1$. (We define $H^1$ as the closure of $H^\infty$ in $L^1=L^1(\T,m)$ and endow it with the $L^1$ norm.) Indeed, consider the subspace $H^1_0:=zH^1$, which is the kernel of the functional 
$$h\mapsto\widehat h(0),\qquad h\in H^1.$$
A well-known theorem of de Leeuw and Rudin \cite{dLR} identifies the extreme points of $\text{\rm ball}(H^1)$ as the unit-norm outer functions in $H^1$. Furthermore, a result of Cima and Thomson \cite{CT} tells us that the strongly extreme points of $\text{\rm ball}(H^1)$ coincide with its extreme points. It follows easily that the extreme, as well as strongly extreme, points of $\text{\rm ball}(H^1_0)$ are of the form $zF$, where $F\in H^1$ is an outer function of norm $1$. Thus, in contrast to the $H^\infty$ situation, the set of extreme---or strongly extreme---points for the smaller space (i.e., $H^1_0$) is disjoint from the corresponding set for the bigger one (i.e., $H^1$). 

\smallskip (4) In the same vein, we mention some other subspaces of $H^1$, more interesting than $H^1_0$ and also defined in spectral terms, whose extreme points have been determined. Namely, given a subset $\La$ of $\Z_+:=\{0,1,2,\dots\}$, we put
$$H^1(\La):=\{f\in H^1:\,\widehat f(k)=0\,\,\text{\rm whenever}\,\,k\notin\La\}.$$
Assuming that $\#(\Z_+\setminus\La)<\infty$, we characterized in \cite{DCRM, DAdv2022} the extreme points of $\text{\rm ball}(H^1(\La))$. Among them, one can always find a non-outer function, so there is no stability result such as Theorem \ref{thm:exstex} in this setting. The geometry of $\text{\rm ball}(H^1(\La))$ was also studied for finite sets $\La$; see \cite{DMRL2000, DAdv2021}.

\smallskip (5) In view of the de Leeuw--Rudin and Cima--Thomson theorems on the (strongly) extreme points of $\text{\rm ball}(H^1)$, the following question arises: Given a unit-norm function $f\in H^1$, can we somehow estimate the quantity 
$$\inf\{\eps>0:\,\eps\,\,\text{\rm is }H^1\text{\rm -admissible for }f\}$$
in terms of the inner factor of $f$? After all, this infimum is zero if and only if $f$ is outer.

\smallskip (6) Finally, it should be mentioned that the stability issue is trivial for the Hardy spaces $H^p$ with $1<p<\infty$. This is due to the uniform convexity of $H^p$ (and in fact of the containing $L^p$ space), which implies that every unit-norm function is a strongly extreme point of the unit ball.

\end{document}